\def\versiondate{18 Oct. 2011}
\input math.macros
\input Ref.macros

\checkdefinedreferencetrue
\continuousnumberingtrue
\continuousfigurenumberingtrue
\theoremcountingtrue
\sectionnumberstrue
\forwardreferencetrue
\citationgenerationtrue
\nobracketcittrue
\hyperstrue
\initialeqmacro

\input\jobname.key
\bibsty{../../texstuff/myapalike}

\def\bfz{{\bf 0}}
\def\Xsp{{\scr X}}  %
\def\Ysp{{\scr Y}}  %
\def\dCov{\mathop {\rm dCov}}
\def\dcov{\mathop {\rm dcov}}
\def\ip#1,#2{\langle #1, #2 \rangle}
\def\Bigip#1,#2{\Big\langle #1, #2 \Big\rangle}
\def\bigip#1,#2{\big\langle #1, #2 \big\rangle}
\def\bry{{\beta}}  %
\def\tr{\mathop {\rm tr}}  %
\def\phih{\widehat\phi}
\def\psih{\widehat\psi}

\def\slightstrut{\vrule height8.5pt depth0pt width0pt}
\def\sphere{\Bbb S}

\def\SRB{\ref b.SRB:measure/, hereinafter referred to as 
{\htmllocref{\bibcode{SRB:measure}}{SRB}}%
\def\SRB{\htmllocref{\bibcode{SRB:measure}}{SRB}}}

\ifproofmode \relax \else\head{}
{Version of \versiondate}\fi 
\vglue20pt

\title{Distance Covariance in Metric Spaces}

\author{Russell Lyons}

\abstract{We extend the theory of distance (Brownian)
covariance from Euclidean
spaces, where it was introduced by Sz\'ekely, Rizzo and Bakirov, to general
metric spaces. We show that for testing independence, it is necessary and
sufficient that the metric space be of strong 
negative type. In particular, we show that this holds for separable Hilbert
spaces, which answers a question of Kosorok. Instead of the manipulations
of Fourier transforms used in the original work, we use elementary
inequalities for metric spaces and embeddings in Hilbert spaces.}

\bottomII{Primary 
62H20,   %
62G20,  	%
51K99.    %
Secondary
 62H15,   	%
30L05.  	%
}
{Negative type, hypothesis testing, independence, distance correlation,
Brownian covariance.}
{Research partially supported by NSF grant DMS-1007244 and Microsoft
Research.}

\bsection{Introduction}{s.intro}

\ref b.SRB:measure/ introduced a new statistical test for the following
problem: given IID samples of a pair of random variables $(X, Y)$, where
$X$ and $Y$ have finite first moments, are $X$
and $Y$ independent?
Among the virtues of their test is that it is extremely simple to compute,
based merely on a quadratic polynomial of the distances between points in
the sample, and that it is consistent against all alternatives (with finite
first moments).
The test statistic is based on a new notion 
called ``distance covariance" or ``distance
correlation".
The paper by \ref b.SR:Brownian/ introduced another new notion, ``Brownian
covariance", and showed it
to be the same as distance covariance. That paper also gave more examples
of its use. This latter paper elicited such interest that it was
accompanied by a 3-page editorial introduction and 42 pages of comments.

Although the theory presented in those papers is very beautiful, it also
gives the impression of being rather technical, relying on various
manipulations with Fourier transforms and arcane integrals.
Answering a question from Sz\'ekely (personal communication, 2010), we show
that almost the entire theory can be developed for general metric spaces,
where it necessarily becomes much more elementary and transparent. A
crucial point of the theory is that the distance covariance of $(X, Y)$ is
0 iff $X$ and $Y$ are independent. This does not hold for general metric
spaces, but we characterize those for which it does hold.
Namely, they are the metric spaces that have what we term ``strong negative
type".

In fact, negative type had arisen already in the work of \SRB.
It was especially prominent in its predecessors,
\refbmulti{SR:hier,SR:normal}.
The notion of strict negative type is standard, but we need a strengthening
of it that we term ``strong negative type". (These notions were conflated
in \SRB\ and \refbmulti{SR:hier,SR:normal}.)

The concept of negative type is old, but has enjoyed a resurgence of
interest recently due to its uses in theoretical computer science, where
embeddings of metric spaces, such as graphs, play a useful role in
algorithms; see, e.g., \ref b.Naor:ICM/ and \ref b.DezaLaurent/.
The fact that Euclidean space has negative type is behind the
following charming and venerable puzzle: 
Given $n$ red points $x_i$ and $n$ blue points $x'_i$ in $\R^p$, show that
the sum $2\sum_{i, j} \|x_i - x'_j\|$ of the distances between the $2n^2$
ordered pairs of points of opposite color is at least the sum $\sum_{i, j}
\big(\|x_i - x_j\| + \|x'_i - x'_j\|\big)$ of the distances between the
$2n^2$ ordered pairs of points of the same color.
The reason the solution is not obvious is that it requires a special
property of Euclidean space.
The connection to embeddings is that, as
\refbmulti{Schoenberg:Ann,Schoenberg:TAMS} showed,
negative type is equivalent to a certain property of embeddability into
Hilbert space. Indeed, if distance in the puzzle were replaced by squared
distance, it would be easy.

If we replace the sums of distances in the puzzle by averages, and then
replace the two finite sets of points by two probability distributions
(with finite first moments), we arrive at an equivalent property, called
negative type. The condition that equality holds only when the two
distributions are equal is called ``strong negative type". It means that a
simple computation involving average distances allows one to distinguish
any two probability distributions. Many statistical tests are aimed at
distinguishing two probability distributions, or distinguishing two
families of distributions. This is what lies directly behind the tests in
\refbmulti{SR:hier,SR:normal}. It is also what lies behind the papers \ref
b.BRS:multi/, \SRB, and \ref b.SR:Brownian/, but there it is somewhat hidden.
We bring this out more clearly in showing how distance covariance allows a
test for independence precisely when the two marginal distributions lie in
metric spaces of strong negative type.

In \ref s.genmetric/, we define distance covariance and prove its basic
properties for general metric spaces. This includes a statistical test for
independence, but it is consistent against all alternatives only in the
case of spaces of strong negative type, as explained in \ref s.negtype/. In
\ref s.negtype/, we also sketch short proofs of Schoenberg's theorem and
short solutions of the above puzzle (none being original).
It turns out that various embeddings
into Hilbert space, though necessarily equivalent at the abstract level,
are useful for different specific purposes.
In both sections, we separate needed results from other interesting results
by putting the latter in explicit remarks.
We show that the full theory extends to separable-Hilbert-space-valued
random variables, which resolves a question of \ref b.Kosorok/.
We remark at the end of the paper that if $(\Xsp, d)$ is a metric space of
negative type, then $(\Xsp, d^r)$ has strong negative type for all $r \in
(0, 1)$; this means that if in a given application one has negative type
but not strong negative type (for example, in an $L^1$ metric space), then
a simple modification of the metric allows the full theory to apply.

\bsection{General Metric Spaces}{s.genmetric}

Let $(\Xsp, d)$ be a metric space.
Let $M(\Xsp)$ denote the finite signed Borel 
measures on $\Xsp$ and $M_1(\Xsp)$ be the
subset of probability measures.
We say that $\mu \in M(\Xsp)$ has a \dfn{finite first moment} if $\int_\Xsp
d(o, x) \,d|\mu|(x) < \infty$ for some $o \in \Xsp$. The choice of
$o \in \Xsp$ does not matter by virtue of the triangle inequality.
If $\mu,
\mu' \in M(\Xsp)$ both have finite first moments, then
$\int d(x, x') \,d\big(|\mu| \times |\mu'|\big)(x, x') < \infty$
since
$d(x, x') \le d(o, x) + d(o, x')$.
Therefore, 
$\int d(x, x') \,d\mu(x) \,d\mu'(x')$
is defined and finite.
In particular, we may define
$$
a_\mu(x) := \int d(x, x') \,d\mu(x')
$$
and
$$
D(\mu) := \int d(x, x') \,d\mu^2(x, x') 
$$
as finite numbers when $\mu \in M(\Xsp)$ has a finite first moment.
Also, write 
$$
d_\mu(x, x') := d(x, x') - a_\mu(x) - a_\mu(x') + D(\mu)
\,.
$$
The function $d_\mu$ is better behaved than $d$ in the following sense:

\procl l.1to2
Let $\Xsp$ be any metric space.
If $\mu \in M_1(\Xsp)$ has a finite first moment, i.e., $d(x, x') \in
L^1(\mu\times\mu)$, then $d_\mu(x, x') \in
L^2(\mu\times\mu)$.
\endprocl

\proof
For simplicity, write $a(x) := a_\mu(x)$ and $a := D(\mu)$.
Let $X, X' \sim \mu$ be independent.
By the triangle inequality, we have
$$
|d(x, x') - a(x)| \le a(x')
\,,
\label e.tri
$$
whence 
$$
\int d_\mu(x, x') \,d\mu^2(x, x') 
=
\Eleft{\big(d(X, X') - a(X) - a(X') + a\big)^2}
\le
\E[X_1 X_2]
\,,
$$
where $X_1 := \max \big\{ |a-2a(X')|, a\big\}$ and $X_2 := \max \big\{
|a-2a(X)|, a\big\}$.
Since $X_1$ and $X_2$ are integrable and independent, $X_1 X_2$ is also
integrable, with $\E[X_1 X_2] \le 4 a^2$. 
\Qed

The proof of \ref l.1to2/ shows that $\|d_\mu\|_2 \le 2 D(\mu) = 2 \|d\|_1$,
but the factor of 2 will be removed in \ref p.CS/.

We call $\mu \in M(\Xsp)$ \dfn{degenerate} if its support consists of only
a single point.

\procl r.dmud
Let $\mu \in M_1(\Xsp)$ have finite first moment and be non-degenerate.
Although $d_\mu(x, x') < d(x, x')$ for all $x, x' \in \Xsp$,
it is not true that
$|d_\mu(x, x')| \le d(x, x')$ for all $x, x'$ in the support of $\mu$.
To see these, we prove first that 
$$
a_\mu(x) > D(\mu)/2
\label e.amu
$$
for all $x \in \Xsp$.
Indeed, $D(\mu) = \int d(x', x'') \,d\mu^2(x', x'') \le \int [d(x', x) + d(x,
x'')] \,d\mu^2(x', x'') = 2 a_\mu(x)$. Furthermore, if equality holds, then
$d(x', x'') = d(x', x) + d(x, x'')$ for all $x', x''$ in the support of $\mu$.
Put $x' = x''$ to get that $x = x'$, contradicting that $\mu$ is not
degenerate. This proves \ref e.amu/.
Using \ref e.amu/ twice in the definition of $d_\mu$ gives $d_\mu < d$.
On the other hand, \ref e.amu/ also shows that
$d_\mu(x, x) < 0 = -d(x, x)$ for all $x$.
\endprocl

Now let $(\Ysp, d)$ be another metric space.
Let $\theta \in M_1(\Xsp \times \Ysp)$ have finite first moments for each of
its marginals $\mu$ on $\Xsp$ and $\nu$ on $\Ysp$.
Define
$$
\delta_\theta\big((x, y), (x', y')\big) := d_\mu(x, x') d_\nu(y, y')
\,.
$$
By \ref l.1to2/ and the Cauchy-Schwarz inequality, we may define
$$
\dcov(\theta)
:=
\int \delta_\theta\big((x, y), (x', y')\big)
\,d\theta^2\big((x, y), (x', y')\big)
\,.
$$
It is immediate from the definition that if $\theta$ is a product measure,
then $\dcov(\theta) = 0$; the converse statement is not always true and is
the key topic of the theory. Metric spaces that satisfy this are
characterized in \ref s.negtype/ as those of strong negative type.
Similarly, spaces for which $\dcov \ge 0$ are characterized in \ref
s.negtype/ as those of negative type.
\SRB\ call the {\it square root\/} of $\dcov(\theta)$ the
\dfn{distance covariance} of $\theta$, but they work only in the context of
Euclidean spaces, where $\dcov \ge 0$. They denote that square root by
$\dCov(\theta)$.

When $(X, Y)$ are random variables with distribution $\theta \in M_1(\Xsp
\times \Ysp)$, we also write $\dcov(X, Y) := \dcov(\theta)$.
If $(X, Y)$ and $(X', Y')$ are independent, both with distribution
$\theta$ having marginals $\mu$ and $\nu$, then
$$
\dcov(\theta)
=
\Ebig{\big(d(X, X') - a_\mu(X) - a_\mu(X') + D(\mu)\big)
\big(d(Y, Y') - a_\nu(Y) - a_\nu(Y') + D(\nu)\big)}
\,.
$$
The following generalizes (2.5) of \SRB. 

\procl p.CS
Let $\Xsp$ and $\Ysp$ be any metric spaces.
Let $\theta \in M_1(\Xsp \times \Ysp)$ have finite first moments for each of
its marginals $\mu$ on $\Xsp$ and $\nu$ on $\Ysp$.
Let $(X, Y) \sim \theta$.
Then
$$\eqaln{
|\dcov(X, Y)| 
&\le
\sqrt{\dcov(X, X)\dcov(Y, Y)}
\label e.CS
\cr&\le D(\mu) D(\nu)
\,.
}$$
Furthermore, $\dcov(X, X) = D(\mu)^2$ iff $\mu$ is concentrated on at most
two points.
\endprocl

\proof
The Cauchy-Schwarz inequality shows \ref e.CS/.
It remains to show that 
$$
\dcov(X, X) \le D(\mu)^2
\label e.1to2
$$
and to analyze the case of
equality. As before, write $a(x)
:= a_\mu(x)$ and $a := D(\mu)$.
By \ref e.tri/, we have 
$$
\Ebig{|d(X, X') - a(X)|a(X)} \le \Ebig{a(X')a(X)}
= a^2 < \infty
\,,
$$
whence $\Ebig{[d(X, X') - a(X)]a(X)} = 0$ by Fubini's theorem (i.e.,
condition on $X$). Similarly,
$\Ebig{[d(X, X') - a(X')]a(X')} = 0$. Thus, expanding the square in
$\dcov(X, X) = \Eleft{\big(d(X, X') - a(X) - a(X') + a\big)^2}$ and
replacing $d(X, X')^2$ there by the larger quantity $d(X, X')[a(X) +
a(X')]$ yields $[d(X, X') - a(X)]a(X) + [d(X, X') - a(X')]a(X')$ plus other
terms that are individually integrable with integrals summing to $a^2$.
This shows the inequality \ref e.1to2/.
Furthermore, it shows that equality holds iff
for all points $x, x'$ in the support of $\mu$, if $d(x, x') \ne 0$, then
$d(x, x') = a(x) + a(x')$. Since the right-hand side equals $\int [d(x, o)
+ d(o, x')] \,d\mu(o)$, it follows that $d(x, x') = d(x, o) + d(o, x')$ for
all $o$ in the support of $\mu$. If there is an $o \ne x, x'$ in the
support of $\mu$, then we
similarly  have that $d(x, o) = d(x, x') + d(x', o)$. Adding these equations
together shows that $d(o, x') = 0$, a contradiction.
That is, if $\dcov(X, X) = D(\mu)^2$, then the support of $\mu$ has size 1
or 2.
The converse is clear.
\Qed

The next proposition generalizes Theorem 4(i) of \SRB.

\procl p.dcovXX0
If $\dcov(X, X) = 0$, then $X$ is degenerate.
\endprocl

\proof
As before, write $a(x) := a_\mu(x)$ and $a := D(\mu)$, where $X \sim \mu$.
The hypothesis implies that
$d(X, X') - a(X) - a(X') + a = 0$ a.s.
Since all functions here are continuous, we have
$d(x, x') - a(x) - a(x') + a = 0$ for all $x, x'$ in the support of $\mu$.
Put $x = x'$ to deduce that for all $x$ in the support of $\mu$, we have
$a(x) = a/2$.
Therefore, $d(X, X') = 0$ a.s. 
\Qed

Assume that $\mu$ and $\nu$ are non-degenerate.
Then the right-hand side of \ref e.CS/ is not 0; the quotient
$\dcov(\theta)/[D(\mu)D(\nu)]$
is the square of what is called the
\dfn{distance correlation} of $\theta$ in \SRB.
In \SRB, this quotient is always non-negative.

This next proposition extends Theorem 3(iii) of \SRB.

\procl p.dcorr1
If $\mu$ and $\nu$ are non-degenerate and
equality holds in \ref e.CS/,
then for some $c > 0$, there is a continuous map $f : \Xsp \to \Ysp$ such
that for all $x, x'$ in the support of $\mu$, we have $d(x, x') = c
d\big(f(x), f(x')\big)$ and $y = f(x)$ for $\theta$-a.e.\ $(x, y)$.
\endprocl

\proof
Write $a(x) := a_\mu(x)$, $a := D(\mu)$, $b(y) := a_\nu(y)$, and $b :=
D(\nu)$.
Equality holds in \ref e.CS/ 
iff there is some constant $c$ such that 
$$
d(x, x') - a(x) - a(x') + a
=
c\big(d(y, y') - b(y) - b(y') + b\big)
$$
for $\theta^2$-a.e.\ $(x, y), (x', y')$,
i.e.,
$$
d(x, x') - c d(y, y') 
=
a(x) - c b(y)
+ a(x') - c b(y') 
+ c b - a
\,.
$$
Since all functions here are continuous, the same holds
for all $(x, y), (x', y')$ in the support of $\theta$.
Put $(x, y) = (x', y')$ to deduce that for all $(x, y)$ in the support of
$\theta$, we have $a(x) - c b(y) = (a - c b)/2$.
This means that $d(x, x') = c d(y, y')$ $\theta^2$-a.s.
The conclusion follows.
\Qed

We now extend Theorem 2 of \SRB.

\procl p.lln
Let $\Xsp$ and $\Ysp$ be metric spaces.
Let $\theta \in M_1(\Xsp \times \Ysp)$ have marginals with finite first
moment.
Let $\theta_n$ be the (random) empirical measure of the first $n$ samples
from an infinite sequence of IID samples of $\theta$.
Then $\dcov(\theta_n) \to \dcov(\theta)$ a.s.
\endprocl

\proof
Let $(X^i, Y^i) \sim \theta$ be independent for $1 \le i \le 6$.
Write 
$$
f(z_1, z_2, z_3, z_4)
:=
d(z_1, z_2) - d(z_1, z_3) - d(z_2, z_4) + d(z_3, z_4)
\,.
$$
Here, $z_i \in \Xsp$ or $z_i \in \Ysp$.
The triangle inequality gives that 
$$
|f(z_1, z_2, z_3, z_4)|
\le
g(z_1, z_3, z_4)
:= \max \big\{ d(z_3, z_4), d(z_1, z_3) \big\}
$$
and
$$
|f(z_1, z_2, z_3, z_4)|
\le
g(z_2, z_4, z_3)
= \max \big\{ d(z_3, z_4), d(z_2, z_4) \big\}
\,.
$$
Since $g(X^1, X^3, X^4)$ and $g(Y^2, Y^6, Y^5)$ are integrable and
independent, it follows that 
$$
h\big((X^1, Y^1), \ldots, (X^6, Y^6)\big)
:=
f(X^1, X^2, X^3, X^4) f(Y^1, Y^2, Y^5, Y^6)
$$
is integrable.
Fubini's theorem thus shows that its expectation equals $\dcov(\theta)$.
Similarly, $\dcov(\theta_n)$ are the V-statistics for the kernel 
$h$ %
of degree 6.
Hence, the result follows. 
\Qed

The proof of \ref p.lln/ for general metric spaces is more straightforward
if second moments are finite, as in Remark 3 of \SRB.

We next extend Theorem 5 of \SRB.

\procl t.cltgen
Let $\Xsp$, $\Ysp$ be metric spaces. 
Let $\theta \in M_1(\Xsp \times \Ysp)$ have 
marginals $\mu$, $\nu$ with finite first moment.
Let $\theta_n$ be the empirical measure of the first $n$ samples from an
infinite sequence of IID samples of $\theta$.
Let $\lambda_i$ be the eigenvalues (with multiplicity)
of the map that sends $F \in L^2(\theta)$ to the function
$$
(x, y) \mapsto \int \delta_\theta\big((x, y), (x', y')\big) F(x', y') 
\,d\theta(x', y')
\,.
$$
If $\theta = \mu\times\nu$, 
then $n \dcov(\theta_n) \Rightarrow \sum_i \lambda_i Z_i^2$,
where $Z_i$ are IID standard normal random variables and $\sum_i
\lambda_i = D(\mu) D(\nu)$.
\endprocl

\proof
We use the same notation as in the proof of \ref p.lln/.
That proof shows that $h$ is integrable when $\mu$
and $\nu$ have finite first moments; the case $X^i = Y^i$ shows then that
$f(X^1, X^2, X^3, X^4)$ has finite second moment. Therefore, when $\theta =
\mu\times\nu$,
$h\big((X^1, Y^1), \ldots, (X^6, Y^6)\big)$ has finite
second moment.

Assume now that $\theta = \mu\times\nu$.
Then kernel $h$ is degenerate of order 1.
Let $\bar h$ be the symmetrized version of $h$. Then since $\theta = \mu
\times \nu$,
$$
\bar h_2\big((x, y), (x', y')\big)
:=
\Ebig{ \bar h\big((x, y), (x', y'), (X^3, Y^3), \ldots, (X^6, Y^6)\big)}
=
\delta_\theta\big((x, y), (x', y')\big)/15
\,.
$$
Hence the result follows from the well-known theory of degenerate
V-statistics (compare Theorem 5.5.2 in \ref b.Serfling/ or Example 12.11 in
\ref b.Vaart/ for the case of U-statistics).
Finally, we have $\sum \lambda_i = \int \delta_\theta\big((x, y), (x, y)\big) 
\,d\theta(x, y) = D(\mu) D(\nu)$ since $\theta =
\mu\times\nu$.
\Qed

\procl c.testgen
Let $\Xsp$, $\Ysp$ be metric spaces.
Let $\theta \in M_1(\Xsp \times \Ysp)$ have non-degenerate
marginals $\mu$, $\nu$ with finite first moment.
Let $\theta_n$ be the empirical measure of the first $n$ samples from an
infinite sequence of IID samples of $\theta$.
Let $\mu_n$, $\nu_n$ be the marginals of $\theta_n$.
If $\theta = \mu \times \nu$, then 
$$
{n \dcov(\theta_n) \over D(\mu_n) D(\nu_n)} \Rightarrow {\sum_i \lambda_i Z_i^2
\over
D(\mu) D(\nu)}
\,,
\label e.testlim
$$
where $\lambda_i$ and $Z_i$ are as in \ref t.cltgen/ and the right-hand
side has expectation 1.
If $\dcov(\theta) \ne 0$,
then the left-hand side of \ref e.testlim/
tends to $\pm\infty$ a.s.
\endprocl

\proof
Since $D(\mu_n)$ and $D(\nu_n)$ are V-statistics, we have
$D(\mu_n) \to D(\mu)$ and $D(\nu_n) \to
D(\nu)$ a.s. Thus, the
first case follows from \ref t.cltgen/. The second case follows
from \ref p.lln/.
\Qed

\procl r.eigs
Since $\theta = \mu \times \nu$, the map in \ref t.cltgen/ is the tensor
product of the maps 
$$
L^2(\mu) \ni F \mapsto \big(x \mapsto \int d_\mu(x, x')
F(x') \,d\mu(x')\big)
$$
and 
$$
L^2(\nu) \ni F \mapsto \big(y \mapsto \int
d_\nu(y, y') F(y') \,d\nu(y')\big)
\,.
$$
Therefore, the eigenvalues $\lambda_i$
are the products of the eigenvalues of these two maps.
\endprocl

\bsection{Spaces of Negative Type}{s.negtype}

\ref c.testgen/ is incomplete in that it does not specify what happens when
$\dcov(\theta) = 0$ and $\theta$ is not a product measure.
In order for the statistics $\dcov(\theta_n)$ to give a test for
independence that is consistent against all alternatives, it suffices to
rule out this missing case. In this section, we show that this case never
arises for metric spaces of strong negative type, but otherwise it does.
This will require the development of several other theorems of independent
interest.
We intersperse these theorems with their specializations to Euclidean
space.

The puzzle we recalled in the introduction can be stated the following
way for a metric space $(\Xsp, d)$: Let
$n \ge 1$ and $x_1, \ldots, x_{2n} \in \Xsp$. Write $\alpha_i$
for the indicator that $x_i$ is red minus the indicator that $x_i$ is blue.
Then $\sum_{i=1}^{2n} \alpha_i = 0$ and
$$
\sum_{i, j \le 2n} \alpha_i \alpha_j d(x_i, x_j) \le 0
\,.
$$
By considering repetitions of $x_i$ and taking limits, we arrive at a
superficially more general property:
For all $n \ge 1$, $x_1, \ldots, x_n \in \Xsp$, and $\alpha_1, \ldots,
\alpha_n
\in \R$ with $\sum_{i=1}^n \alpha_i = 0$, we have 
$$
\sum_{i, j \le n} \alpha_i \alpha_j d(x_i, x_j) \le 0
\,.
\label e.ntdef
$$
We say that $(\Xsp, d)$ has \dfn{negative type}
if this property holds.
A list of metric spaces of negative type appears as Theorem 3.6 of \ref
b.Meckes:pdms/; in particular, this includes
all $L^p$ spaces for $1 \le p \le 2$.
On the other hand, $\R^n$ with the $\ell^p$-metric
is not of negative type whenever $3 \le n \le \infty$
and $2 < p \le \infty$, as proved by 
\ref b.Dor/ combined with Theorem 2 of \ref b.BDK/;
see \ref b.KolLon/ for an extension to spaces that include some Orlicz
spaces, among others.

If we define the $n \times n$ matrix $K$ whose $(i, j)$ entry is $d(x_i,
x_j)$, then \ref e.ntdef/ says, by definition, that $K$ is conditionally
negative semidefinite. 
This explains the name ``negative type".
We can construct another matrix $\bar K$ from $K$ that is negative
semidefinite as follows: Let $P$ be the orthogonal projection of
$\R^n$ onto the orthocomplement of the constant vectors. Then as operators,
$\bar K := P K P$.
Let $\mu_n$ be the empirical measure of $x_1, \ldots, x_n$.
The $(i, j)$ entry of $\bar K$ is easily verified to be
$d_{\mu_n}(x_i, x_j)$, which begins to
explain the appearance of $d_\mu$ in \ref s.genmetric/.
We write $\bar K \le 0$ to mean that $\bar K$ is negative semidefinite.

If $\Xsp$ and $\Ysp$ are both metric spaces of negative type and $(x_i,
y_i) \in \Xsp \times \Ysp$, then
let $K$ and $L$ be the distance matrices for $x_i$ and $y_i$, respectively.
Let $\theta_n$ be the empirical measure of the sequence
$\Seq{(x_i, y_i) \st 1 \le i \le n}$.
We have $\bar K \le 0$ and $\bar L \le 0$, whence $\tr \big(\bar K \bar
L\big) = \tr
\Big(
\sqrt{\slightstrut -\bar K} \sqrt{\slightstrut -\bar L}
\sqrt{\slightstrut -\bar L} \sqrt{\slightstrut -\bar K} 
\Big)
\ge 0$.
That is,
$$
0 \le
\tr \big(\bar K \bar L \big)
=
n^2 \dcov(\theta_n)
\,.
$$
This begins to explain the origin of $\dcov$.
To go further, we use embeddings into Hilbert space.

Now $\Xsp$ is of negative type iff there is a Hilbert space $H$ and a
map $\phi : \Xsp \to H$ such that $\all {x, x' \in \Xsp} d(x, x') =
\|\phi(x) - \phi(x')\|^2$, as shown by
\refbmulti{Schoenberg:Ann,Schoenberg:TAMS}.
Indeed, given such a $\phi$, \ref e.ntdef/ is easy to verify: see \ref
e.negbar/ below. For the converse, consider $x_1, \ldots, x_n \in \Xsp$.
Since $\bar K \le 0$, there are vectors $v_i \in \R^n$ such that $\ip v_i,
{v_j}$ is the $(i, j)$-entry of $-\bar K$ for all $i, j$ (the matrix
$\sqrt{\slightstrut -\bar K}$ has $v_i$ for its $i$th column). Computing
$\|v_i - v_j\|^2$ then yields $\|v_i/\sqrt2 - v_j/\sqrt2\|^2 = d(x_i, x_j)$.
This provides a map $\phi$ defined on the points $x_1, \ldots, x_n$.
When we increase the domain of such a $\phi$, the distances of the images
already defined are preserved, whence we may embed all these images in a
fixed Hilbert space. 
If $\Xsp$ is separable, we may thus define $\phi$ on a countable
dense subset by induction, and then extend by continuity.
In general and alternatively, define 
$$
d_o(x, x') := \big[d(x, o) + d(o, x') - d(x, x')\big]/2
$$
for some fixed $o \in \Xsp$.
Let $V$ be the finitely supported functions on $\Xsp$.
The fact that $\Xsp$ is of negative type implies that $\ip f, g := \sum_{x,
x' \in \Xsp} f(x) g(x') d_o(x, x')$ is a semi-inner product on $V$.
The Cauchy-Schwarz inequality implies that
$V_0 := \big\{f \in V \st \ip f, f = 0\big\}$ is a subspace of $V$. Let $H$
be the completion of $V/V_0$. Then the map $\phi : x \mapsto \II x + V_0$
has the property desired.
Note that $H$ is separable when $\Xsp$ is.

Of course, any two isometric embeddings $\phi_1, \phi_2 : (\Xsp, d^{1/2})
\to H$ are equivalent in the sense that there exists an isometry $g : H_1 \to
H_2$ such that $\phi_2 = g \circ \phi_1$, where $H_i$ is the closed affine
span of the image of $\phi_i$: Define $g\big(\phi_1(x)\big) :=
\phi_2(x)$ for $x \in \Xsp$, extend by affine
linearity (which is well defined by 
a property of Euclidean space), and then extend by continuity.
We shall call
an isometric embedding $\phi : (\Xsp, d^{1/2}) \to H$ simply an
\dfn{embedding}.

A direct proof that $\R^n$ is of negative type is the following. When $n =
1$, define $\phi(x)$ to be the function 
$\I {[0, \infty)} - \I {[x, \infty)}$ in $L^2(\R, \lambda)$,
where $\lambda$ is Lebesgue measure. This is
easily seen to have the desired property.
When $n \ge 2$,
define $f_x(s) := \|x - s\|^{-(n-1)/2}$ and $g_x := f_x - f_\bfz$
for $x \in \R^n$. 
Then $g_x \in L^2(\R^n, \lambda^n)$, as
calculus shows (for large $s$, we have $g_x(s) =
O\big(\|s\|^{-(n+1)/2}\big)$).
Furthermore, there is a constant $c$ such that
$\|g_x\|_2 = c\|x\|^{1/2}$ by homogeneity, whence translation
invariance gives $\|g_x - g_{x'}\|_2 = \|g_{x - x'}\|_2 =
c \|x - x'\|^{1/2}$, so that $\phi(x) := g_x/c$ has the desired property.
Call this embedding the \dfn{Riesz embedding} since $f_x(s)$ is a Riesz
kernel.

Another embedding $\phi$ for $\R^n$ is as follows: $\phi(x)$ is the
function $s \mapsto c(1 - e^{- i s \cdot x})$ in $L^2(F \lambda^n)$ for some
constant $c$, where $F(s) := \|s\|^{-(n+1)}$.
See Lemma 1 of \ref b.SR:hier/ for a proof.
This is the Fourier transform of the Riesz embedding, in other words, the
composition of the Riesz embedding with the Fourier isometry.
We shall refer to this embedding as the \dfn{Fourier embedding}.

Other important embeddings use Brownian motion.
When $n=1$,
let $B_x$ be Brownian motion defined for $x \in \R$ with
$B_0 = \bfz$. We may then define
$\phi(x) := B_x$, thought of as a function in $L^2(\P)$
for some probability measure $\P$. Likewise, the case $n \ge 2$ can be
accomplished by using L\'evy's multiparameter Brownian motion.
We shall refer to these embeddings as the \dfn{Brownian embeddings}.
Sample-path continuity of these Brownian motions plays no role
for us; only their Gaussian structure matters.
In fact, their existence depends only on the fact that $\R^n$ has negative
type.

An embedding that does not rely on calculation goes as follows: Let
$\sigma$ be the (infinite) Borel measure on half-spaces $S \subset\R^n$
that is invariant under translations and rotations, normalized so that
$$
\sigma\big(\{\bfz \in S, x \notin S\}\big) = \|x\|/2
\label e.cut
$$
for $\|x\|= 1$.
If we parametrize half-spaces as $S = \{ x \in \R^n \st z \cdot x \le s\}$
with $z \in \sphere^{n-1}$ and $s \in \R$, then  $\sigma = c_n \Omega_n
\times \lambda$ for some constant $c_n$, where $\Omega_n$ is volume measure
on $\sphere^{n-1}$.
Scaling shows that \ref e.cut/ holds for all $x$.
Now let $\phi(x)$ be the function $S \mapsto \I S(\bfz) - \I S(x)$ 
in $L^2(\sigma)$. We call this the \dfn{Crofton embedding}, as \ref
b.Crofton/ was the first to give a formula for the distance of points in
the plane in terms of lines intersecting the segment joining them.

We return now to general metric spaces of negative type.
Suppose that $\mu_1, \mu_2 \in M_1(\Xsp)$ have finite first
moments.
By approximating $\mu_i$ by probability measures of finite support (e.g.,
IID samples give V-statistics), we see
that when $\Xsp$ has negative type,
$$
D(\mu_1 - \mu_2) \le 0
\,.
\label e.negtype
$$
We say that $(\Xsp, d)$ has \dfn{strong negative type} if it has negative
type and equality holds in \ref e.negtype/ only when $\mu_1 = \mu_2$. When
$\mu_i$ are restricted to measures of finite support, then this is the
condition that $(\Xsp, d)$ be of \dfn{strict negative type}. 
A simple example of a metric space of non-strict negative type is $\ell^1$
on a 2-point space, i.e., $\R^2$ with the $\ell^1$-metric.

Consider an embedding $\phi$ as above.
Define the (linear) barycenter
map $\bry = \bry_\phi: \mu \mapsto \int \phi(x) \,d\mu(x)$
on the set of measures $\mu \in M(\Xsp)$ with finite first moment.
(Although it suffices that $\int d(o, x)^{1/2} \,d|\mu|(x)  < \infty$ to
define $\bry(\mu)$, this will not suffice for our purposes.)
Note that
$$
\int\!\!\int d(x_1, x_2) \,d\mu_1(x_1) \,d\mu_2(x_2) 
=
-2 \bigip \bry(\mu_1), {\bry(\mu_2)}
$$
when $\mu_i \in M(\Xsp)$ satisfy $\mu_i(\Xsp) = 0$.
In particular,
$$
D(\mu)
=
-2 \| \bry(\mu) \|^2
\label e.negbar
$$
when $\mu \in M(\Xsp)$ satisfies $\mu(\Xsp) = 0$.
Thus:

\procl p.strnegbar
Let $\Xsp$ have negative type as witnessed by the embedding $\phi$.
Then $\Xsp$ is of strong negative type iff the barycenter map $\bry_\phi$
is injective on the set of probability measures on $\Xsp$ with finite first
moment.
\endprocl

For example, Euclidean spaces have strong negative type; this is most
directly seen via the Fourier embedding, since then $\bry(\mu)$ is the
function $s \mapsto c\big(1- \widehat\mu(s)\big)$, where $\widehat\mu$ is
the Fourier transform of $\mu \in M_1(\R^n)$.
Alternatively, one can see this via the Crofton embedding and the
Cram\'er-Wold device, but the only decent proof of that device uses Fourier
transforms. (Of course, in one dimension, the Crofton embedding is
simple and easily shows that $\R$ has strong negative type without the use
of Fourier transforms.)
The barycenter of $\mu$ for the Riesz embedding is essentially the Riesz
potential of $\mu$; more precisely, if $\mu$ and $\mu'$ are probability
measures with finite first moment, then up to a constant factor,
$\bry(\mu-\mu')$ is the Riesz potential of $\mu - \mu'$ for the exponent
$(n-1)/2$.

\procl r.M1metric
Another way of saying \ref p.strnegbar/ is that a metric space $(\Xsp, d)$
has strong negative type iff the map $(\mu_1, \mu_2) \mapsto \sqrt{-D(\mu_1
- \mu_2)/2}$ is a metric on the set of probability measures on $\Xsp$ with
finite first moment, in which case it extends the metric on $(\Xsp,
d^{1/2})$ when we identify $x \in \Xsp$ with the point mass at $x$.
\endprocl

\procl r.negstrong
Here we give an example of a metric space of strict negative type that is
not of strong negative type. In fact, it fails the condition for
probability measures with countable support.
The question amounts to
whether, given a subset of a Hilbert
space in which no 3 points form an obtuse triangle and such that the
barycenter of every finitely supported probability measure determines the
measure uniquely, the barycenter of every probability measure determines the
measure uniquely. The
answer is no. For example, let $\Seq{e_i}$ be an orthonormal basis of
a Hilbert space. The desired subset consists of the vectors
$$\displaylines{
e_1 \,,\cr
e_1 + e_2/2 \,,\cr
e_2 + e_3 \,,\cr
e_3 + e_4/2 \,,\cr
e_4 + e_5 \,,\cr
e_5 + e_6/2 \,,\cr
}$$
etc.
It is obvious that finite convex combinations are unique and 
that there are no obtuse angles. But if $v_n$ denotes the 
$n$th vector, then
$$
v_1/2 + v_3/4 + v_5/8 + \cdots
=
v_2/2 + v_4/4 + v_6/8 + \cdots
\,.
$$
\endprocl

\procl r.distancetransform
If $\Xsp$ is a metric space of negative type, then $\alpha :
\mu \mapsto a_\mu$ is
injective on $\mu \in M_1(\Xsp)$ with finite first moment
iff $\Xsp$ has strong negative type.
Part of this statement is contained in Theorem 3.6 of \ref b.NickWolf/.
To prove it, let $\phi$ be an embedding of $\Xsp$ such that $\bfz$ lies in the
image of $\phi$, which we may achieve by translation. Then 
$$
a_\mu(x) = \|\phi(x)\|^2 - 2
\bigip \phi(x), {\bry(\mu)} + \int \|\phi(x')\|^2 \,d\mu(x')
\,,
$$
whence
$a_\mu = a_{\mu'}$ iff $\bigip \phi(x), {\bry(\mu)} = \bigip \phi(x),
{\bry(\mu')}$ for all $x$ (first use $x$ so that $\phi(x) = \bfz$) iff 
$\bigip z, {\bry(\mu)} = \bigip z, {\bry(\mu')}$ for all $z$ in the closed
linear span of the image of $\phi$ iff $\bry(\mu) = \bry(\mu')$. Now apply
\ref p.strnegbar/.
On the other hand, there are metric spaces not of negative type for which
$\alpha$ is injective on the probability measures: e.g., take a finite
metric space in which the distances to a fixed point are linearly
independent.
The map $\alpha$ is injective also for all separable $L^p$ spaces ($1 < p <
\infty$): see \ref b.Linde:uniq/ or \ref b.GK:pot/.
\endprocl

Given an $H$-valued random variable $Z$ with finite first moment, we
define its \dfn{variance} to be $\Var(Z) := \Ebig{\|Z - \E[Z]\|^2}$.

\procl p.varX
If $\Xsp$ has negative type as witnessed by the embedding $\phi$ and $\mu
\in M_1(\Xsp)$ has finite first moment, then for all $x, x' \in \Xsp$,
$$
a_\mu(x) = \|\phi(x) - \bry_\phi(\mu)\|^2 + D(\mu)/2
\,,
$$
$D(\mu) = 2 \Var\big(\phi(X)\big)$ if $X \sim \mu$, and 
$$
d_\mu(x, x') = -2 \bigip {\phi(x) - \bry_\phi(\mu)},
{\phi(x') - \bry_\phi(\mu)}
\,.
$$
\endprocl

\proof
Let $X \sim \mu$.
We have 
$$\eqaln{
a_\mu(x) 
&=
\Ebig{d(x, X)}
=
\Ebig{\|\phi(x) - \phi(X)\|^2}
=
\EBig{\big\|\big(\phi(x) - \bry(\mu)\big) - \big(\phi(X) -
\bry(\mu)\big)\big\|^2}
\cr&=
\|\phi(x) - \bry_\phi(\mu)\|^2 + 
\Var\big(\phi(X)\big)
\,.
}$$
Integrating over $x$ gives the first two identities. Substituting the first
identity into the definition of $d_\mu$ gives the last identity. 
\Qed

For simplicity, we may, without loss of generality,
work only with real Hilbert spaces.
Let $\Xsp$ and $\Ysp$ be metric spaces of negative type, witnessed by the
embeddings $\phi$ and $\psi$, respectively.
Consider the tensor embedding $(x, y) \mapsto \phi(x) \otimes \psi(y)$ of
$\Xsp \times \Ysp \to H \otimes H$.
This will be the key to analyzing when $\dcov(\theta) = 0$.
Recall that the inner product on $H \otimes H$ satisfies $\bigip {h_1
\otimes h'_1},
{h_2 \otimes h'_2} := \ip h_1, {h_2} \ip h'_1, {h'_2}$.

\procl r.prodmetric
Although we shall not need it, we may
give $\Xsp \times \Ysp$ the associated ``metric"
$$
d_{\phi \otimes \psi}\big((x, y), (x', y')\big)
:=
\|\phi(x) \otimes \psi(y) - \phi(x') \otimes \psi(y')\|^2
\,,
$$
so necessarily it is of negative type.
Actually, one can check that this need not satisfy the triangle inequality
when the origin is not in the images of $\phi$ and $\psi$, but, following a
suggestion of ours, Leonard
Schulman (personal communication, 2010)
showed that it is indeed a metric when the images of $\phi$ and
$\psi$ both contain the origin.
Since we may translate $\phi$
and $\psi$ so that this holds, we may take this to be a metric if we wish.
In this case, one can also express $d_{\phi \otimes \psi}$ in terms of the
original metrics on $\Xsp$ and $\Ysp$.
\endprocl

\procl p.tendcov 
Let $\Xsp$, $\Ysp$ have negative type as witnessed by the embeddings
$\phi$, $\psi$.
Let $\theta \in M_1(\Xsp \times \Ysp)$ have marginals
$\mu \in M_1(\Xsp)$ and $\nu \in M_1(\Ysp)$, both with finite first moment.
Then $\theta \circ (\phi \otimes \psi)^{-1}$ has finite first moment, so
that $\bry_{\phi \otimes \psi}(\theta)$ is defined, and
we have that 
$$
\dcov(\theta) = 4\big\|\bry_{\phi \otimes \psi}(\theta-\mu \times \nu)\big\|^2
\,.
$$
\endprocl

\proof
Write $\phih := \phi - \bry_\phi(\mu)$ and
$\psih := \psi - \bry_\psi(\nu)$.
By \ref p.varX/, we have
$$\eqaln{
\dcov(\theta)
&=
4 \int \bigip \phih(x), {\phih(x')} \bigip \psih(y), {\psih(y')}
\,d\theta^2\big((x, y), (x', y')\big)
\cr&=
4 \int \bigip {\phih(x) \otimes \psih(y)}, {\phih(x') \otimes \psih(y')}
\,d\theta^2\big((x, y), (x', y')\big)
\cr&=
4 \|\bry_{\phih \otimes \psih}(\theta)\|^2
\,.
}$$
In addition, since $\|\phi(x)\| \in L^2(\mu)$ and $\|\psi(y)\| \in
L^2(\nu)$, we have $\|\phi(x) \otimes \psi(y)\| \in L^1(\theta)$ by the
Cauchy-Schwarz inequality, whence
$\bry_{\phi \otimes \psi}(\theta)$ is defined and
$$\eqaln{
\bry_{\phih \otimes \psih}(\theta)
&=
\int \phih(x) \otimes \psih(y) \,d\theta(x, y)
=
\int \big(\phi(x) - \bry_\phi(\mu)\big) \otimes
\big(\psi(y) - \bry_\psi(\nu)\big) \,d\theta(x, y)
\cr&=
\int \phi(x) \otimes \psi(y) \,d\theta(x, y) - \bry_\phi(\mu) \otimes
\bry_\psi(\nu)
=
\bry_{\phi \otimes \psi}(\theta-\mu \times \nu)
\,.
\Qed
}$$

In the special case where $\Xsp$ and $\Ysp$ are Euclidean spaces and the
embeddings $\phi, \psi$ are the Fourier embeddings, \ref
p.tendcov/ shows that $\dcov$ coincides with
(the square of) the 
original definition of distance covariance in \SRB\ (see
(2.6) there); while if the embeddings are
the Brownian embeddings, then \ref p.tendcov/ shows
that distance covariance is the same as Brownian covariance (Theorem 8 of
\ref b.SR:Brownian/; the condition there that $X$ and $Y$ have finite
second moments is thus seen to be superfluous).
The Crofton embedding gives 
$$
\bry_{\phi \otimes \psi}(\theta - \mu \times \nu) : 
(z, s, w, t)
\mapsto
c_p c_q
\Big[\theta\big(z \cdot x \le s , w \cdot y \le t \big) - \mu\big(z \cdot x \le s
\big) \nu\big(w \cdot y \le t \big)\Big]
$$
for $\theta \in M_1(\R^p \times \R^q)$ with marginals $\mu, \nu$
having finite first moments,
whence for $(X, Y) \in \R^p \times \R^q$, \ref p.tendcov/ shows that
$$\displaylines{
\dcov(X, Y)
= 
\hfill\break\cr
\hfill
4 c_p c_q
\int \!\!\int
\big|\P[z \cdot X \le s , w \cdot Y \le t ] - \P[z \cdot X \le s ] \P[w \cdot
Y \le t ]\big|^2 \,d(\Omega_p \times \Omega_q) (z, w) \,d\lambda^2(s, t)
\,.
\cr}
$$
When $p = q = 1$, this formula was shown to us by G\'abor Sz\'ekely
(personal communication, 2010).

Write $M^1(\Xsp)$ for the subset of $\mu \in M(\Xsp)$ such that $|\mu|$ has
a finite first moment.
Write $M^{1, 1}(\Xsp \times \Ysp)$ for the subset of $\theta \in M(\Xsp
\times \Ysp)$ such that both marginals of $|\theta|$ have finite first
moment.

\procl l.tenneg
Let $\Xsp$, $\Ysp$ have negative type as witnessed by the embeddings
$\phi$, $\psi$.
If $\phi$ and $\psi$ have the property that
$\bry_\phi$ and $\bry_\psi$ are 
injective on both $M^1(\Xsp)$ and $M^1(\Ysp)$ (not merely on the probability
measures), then $\bry_{\phi \otimes \psi}$ is injective on $M^{1, 1}(\Xsp
\times \Ysp)$.
\endprocl

\proof
Let $\theta \in M^{1, 1}(\Xsp \times \Ysp)$ satisfy
$\bry_{\phi \otimes \psi}(\theta) = 0$. For $k \in H$,
define the bounded linear map $T_k : H \otimes H \to H$ by 
linearity, continuity, and
$$
T_k(u \otimes v) := \ip u, k v
\,.
$$
More precisely, one uses the above definition on $e_i \otimes e_j$ for an
orthonormal basis $\{e_i\}$ of $H$ and then extends.
Also, define
$$
\nu_k(B) := \int \ip \phi(x), k \I B(y) \,d\theta(x, y)
\qquad (B \subseteq \Ysp \hbox{ Borel})
\,,
$$
so that 
$$
\bry_\psi(\nu_k) 
= 
\int \ip \phi(x), k \psi(y) \,d\theta(x, y)
=
\int T_k\big(\phi(x) \otimes \psi(y)\big) \,d\theta(x, y)
=
T_k\big(\bry_{\phi \otimes \psi}(\theta)\big)
= 0
\,.
$$
This implies that $\nu_k = 0$ by injectivity of $\bry_\psi$.
As this is valid for each $k \in H$, we obtain that for every Borel $B
\subseteq \Ysp$, 
$$
\int \phi(x) \I B(y) \,d\theta(x, y)
= 0
\,.
$$
Defining 
$$
\mu_B(A) := \theta(A \times B)
\qquad (A \subseteq \Xsp \hbox{ Borel})
\,,
$$
we have $\bry_\phi(\mu_B) = \int \phi(x) \I B(y) \,d\theta(x, y) =
0$, whence $\mu_B = 0$ by injectivity of $\bry_\phi$.
In other words $\theta(A \times B) = 0$ for every
pair of Borel sets $A$ and $B$. Since such product sets generate the product
$\sigma$-field on $\Xsp \times \Ysp$, it follows that $\theta = 0$.
\Qed

\procl l.tenneg2
Let $\Xsp$ have strong negative type. There exists an embedding
$\phi$ so that 
$\bry_\phi$ is 
injective on $M^1(\Xsp)$ (not merely on the probability
measures).
\endprocl

\proof
If $\phi : \Xsp \to H$ is an embedding that induces
an injective barycenter map on $M^1_1(\Xsp)$, then the map $x \mapsto
\big(\phi(x), 1\big)
\in H \times \R$ is an embedding that induces an
injective barycenter map on $M^1(\Xsp)$.
\Qed

\procl r.prodstrong
We may choose the embeddings so that $d_{\phi \otimes \psi}$ is a metric and
$\bry_{\phi \otimes \psi}$ is injective on $M^1(\Xsp \times
\Ysp)$, which yields that $d_{\phi \otimes \psi}$ is of strong negative
type by \ref p.strnegbar/.
Indeed, first translate $\phi$ and $\psi$ so that each contains $\bfz$ in
its image. This makes $d_{\phi \otimes \psi}$ a metric by \ref
r.prodmetric/. Then
use the embedding $x \mapsto \big(\phi(x), 1\big)$ and likewise for
$\psi$. This does not change the metric.
\endprocl

As we observed in \ref s.genmetric/,
it is immediate from the definition that if $\theta$ is a product measure,
then $\dcov(\theta) = 0$. A converse and the key result of the theory holds
for metric spaces of strong negative type:

\procl t.dcov Suppose that both $\Xsp$ and $\Ysp$ have strong negative
type and $\theta$ is a probability measure on $\Xsp \times \Ysp$ whose
marginals have finite first moment.
If $\dcov(\theta) = 0$, then $\theta$ is a product measure.
\endprocl

This is an immediate corollary of \ref p.tendcov/ and Lemmas \briefref 
l.tenneg/ and \briefref l.tenneg2/.
Therefore, \ref c.testgen/ gives a test for independence that is consistent
against all alternatives when $\Xsp$ and $\Ysp$ both have strong negative
type.
See Theorem 6 of \SRB\ for the significance levels of the
test.

For the Fourier embedding of Euclidean space, \ref t.dcov/
amounts to the fact that
$\theta = \mu \times \nu$ if the Fourier transform of $\theta$ is the
(tensor) product of the Fourier transforms of $\mu$ and $\nu$. This was the
motivation presented in \SRB\ for $\dCov$.

\procl r.categorical
In the case of categorical data, we may embed each data space as a simplex
with edges of unit length. Let the corresponding Hilbert-space vectors be
$e_x/\sqrt 2$ and $f_y/\sqrt 2$, where $e_x$ are orthonormal and $f_y$ are
orthonormal. The product space then embeds as a simplex on the orthogonal
vectors $e_x \otimes f_y/2$ and the barycenter of $\theta$ is $\sum_{x, y}
\theta(x, y) e_x \otimes f_y/2$.
Let $\theta_n$, $\mu_n$, and $\nu_n$ be the empirical
measures as in \ref c.testgen/. 
\ref p.tendcov/ yields
$$
{\dcov}(\theta_n)
=
\sum_{x, y} \big[\theta_n(x, y) - \mu_n(x)\nu_n(y)\big]^2
\,.
$$
The test statistic in \ref e.testlim/ is thus 
$$
{n \sum_{x, y} \big[\theta_n(x, y) - \mu_n(x)\nu_n(y)\big]^2 \over
\sum_x \mu_n(x)[1-\mu_n(x)]
\sum_y \nu_n(y)[1-\nu_n(y)]}
\,.
$$
For comparison, Pearson's $\chi^2$-statistic is
$$
n \sum_{x, y} {\big[\theta_n(x, y) - \mu_n(x)\nu_n(y)\big]^2
\over \mu_n(x)\nu_n(y)}
\,.
$$
\endprocl

\procl r.covdist
As G\'abor Sz\'ekely has remarked (personal communication, 2010), there is
a 2-dimensional random variable $(X, Y)$ such that $X$ and $Y$ are not
independent, yet if $(X', Y')$ is an independent copy of $(X, Y)$, then 
$|X - X'|$ and $|Y - Y'|$ are uncorrelated. Indeed, consider the density
function $p(x, y) := \big(1/4 - q(x) q(y)\big) \I {[-1, 1]^2}(x, y)$ with
$q(x) := -(c/2) \I {[-1, 0]} + (1/2) \I {(0, c)}$, where $c := \sqrt 2 -
1$. Then it is not hard to check that this gives such an example.
\endprocl

\procl r.Direct
According to \ref p.tendcov/, $\dcov(\theta) = -2 D(\theta - \mu \times
\nu)$ for the metric space $(\Xsp \times \Ysp, d_{\phi \otimes \psi})$.
Since this metric space has strong negative type when $\Xsp$ and $\Ysp$ do,
we can view the fact that $\dcov(\theta) = 0$ only for product measures as
a special case of the fact that $D(\theta_1 - \theta_2) = 0$ only when
$\theta_1 = \theta_2$ for $\theta_i \in M^1_1(\Xsp \times \Ysp)$.
Similarly, any other metric on $\Xsp \times \Ysp$ of strong negative type
could be used to give a test of independence via $D(\theta - \mu \times
\nu)$; indeed, when $\Xsp = \R^p$ and $\Ysp = \R^q$, the Euclidean metric
on $\R^{p+q}$ was used by \ref b.BRS:multi/ for precisely such a test.
\endprocl

No such result as \ref t.dcov/ holds if either $\Xsp$ or $\Ysp$ is not
of strong negative type:

\procl p.necnegtype
If $\Xsp$ is not of negative type, then for every metric space
$\Ysp$ with at least two points, there exists $\theta \in
M_1(\Xsp\times\Ysp)$ whose marginals have finite first moments and such
that $\dcov(\theta) < 0$.
If $\Xsp$ is not of strong negative type, then for every metric space
$\Ysp$ with at least two points, there exists $\theta \in
M_1(\Xsp\times\Ysp)$ whose marginals have finite first moments and such
that $\dcov(\theta) = 0$, yet $\theta$ is not a product measure.
\endprocl

\proof
Choose two distinct points $y_1, y_2 \in \Ysp$.
Let $\mu_1 \ne \mu_2 \in M_1(\Xsp)$ have finite first moments and satisfy
$D(\mu_1 - \mu_2) \ge 0$, where $> 0$ applies if $\Xsp$ does not have
negative type.
In this latter case, 
set $\theta := \big(\mu_1 \times \delta(y_1) + \mu_2 \times
\delta(y_2)\big)/2$.
Then a little algebra reveals that
$$
\dcov(\theta)
=
- d(y_1, y_2) D(\mu_1 - \mu_2)/8
<0
\,.
$$

In general, note that
if $x_1 \ne x_2$, then $D\big(\delta(x_1) - \delta(x_2)\big) < 0$, whence
there is some $\gamma \in (0, 1]$ such that if $\tau_i := \gamma \mu_i + (1
- \gamma) \delta(x_i)$, then
$D(\tau_1 - \tau_2) = 0$. 
Set $\theta := \big(\tau_1 \times \delta(y_1) + \tau_2 \times
\delta(y_2)\big)/2$.
Then 
$$
\dcov(\theta)
=
- d(y_1, y_2) D(\tau_1 - \tau_2)/8
=0
\,,
$$
yet $\theta$ is not a product measure.
\Qed

There remains the possibility that the kernel $h$ in the proof of \ref
p.lln/ is degenerate of order 1 only when $\theta$ is a product measure. If
that is true, then \ref c.testgen/ gives a consistent test for independence
even in metric spaces not of negative type, since when $h$ is not
degenerate and $\dcov(\theta) = 0$, $\sqrt n\dcov(\theta_n)$ has a
non-trivial normal
limit in distribution, whence $n\dcov(\theta_n) \to \infty$ a.s. We have
not investigated this possibility.

Since every Euclidean space is of strict negative type, so is every Hilbert
space. Separable Hilbert spaces are even of strong
negative type, though this is considerably more subtle.
Therefore, $\dcov(\theta) = 0$ implies that $\theta \in M_1(\Xsp \times
\Ysp)$ is a product measure when $\Xsp$ and $\Ysp$ are separable Hilbert
spaces, which resolves a question of \ref b.Kosorok/.

\procl t.Hilbertstrong
Every separable Hilbert space is of strong negative type.
\endprocl

\proof
This follows from \ref r.distancetransform/ and Theorem 6 of \ref
b.Linde:rudin/ or Theorem 1 of \ref b.Koldob82/, who prove more.
Likewise, separable $L^p$ spaces with $1 < p < 2$ are of strong negative
type.
However, we give a direct proof that is shorter, which keeps our paper
self-contained.

Our proof relies on a known Gaussian variant of the Crofton embedding.
Let $Z_n$ ($n \ge 1$) be IID standard normal random
variables with law $\rho$ on
$\R^\infty$. Given $u = \Seq{u_n \st n \in
\Z^+} \in \ell^2(\Z^+)$, define the random variable $Z(u) := \sum_{n \ge 1}
u_n Z_n$. Then $Z(u)$ is a centered normal random variable with standard
deviation equal to $\|u\|_2$. Therefore, $\Ebig{|Z(u)|} = c \|u\|_2$ with
$c := \Ebig{|Z_1|}$.

Let $\lambda$ be Lebesgue measure on $\R$.
For $w, u \in \R^\infty$, write $w(u) := \limsup_N \sum_{n=1}^N u_n
w_n$.
We choose $\ell^2(\Z^+)$ as our separable Hilbert space, which we
embed into another Hilbert space,
$L^2(\R^\infty \times \R, \rho \times \lambda)$, by
$$
\phi(u) : (w, s) \mapsto \I {[w(u)/c, \infty)}(s) - \I {[0, \infty)}(s) 
\,.
$$
Then $\|\phi(u) - \phi(u')\|_2^2 = 
\|\phi(u) - \phi(u')\|_1 = \|u - u'\|_2$ for all $u, u' \in
\ell^2(\Z^+)$.
Let $\mu_1, \mu_2 \in M_1\big(\ell^2(\Z^+)\big)$ have finite first
moments. Set $\mu := \mu_1 - \mu_2$.
Because $\int d\mu = 0$, we have
$$
\bry_\phi(\mu) : (w, s) \mapsto \mu\{ u \st w(u) \le c s\}
\,.
$$
Note that since for every $u \in \ell^2(\Z^+)$, the series $w(u)$ converges
$\rho$-a.s., Fubini's theorem tells us that for $\rho$-a.e.\ $w$, $w(u)$
converges for $\mu_i$-a.e.\ $u$.
We need to show that if $\bry_\phi(\mu) = 0$ $\rho \times \lambda$-a.s.,
then $\mu = 0$.
So assume that $\bry_\phi(\mu) = 0$ $\rho \times \lambda$-a.s.
It suffices to show that $\mu\{ u \st \ip u, v \le s\} = 0$ for every
finitely supported $v \in \R^\infty$ and every $s \in \R$, since that
implies that the finite dimensional marginals of $\mu$ are 0 by the
Cram\'er-Wold device.

Let $K \ge 1$. For $w \in \R^\infty$, write $w_{\le K}$ for the vector
$(w_1, \ldots,
w_K) \in \R^K$ and $w_{>K}$ for $(w_{K+1}, w_{K+2}, \ldots ) \in
\R^\infty$. Since the law $\rho$
of $w = \big(w_{\le K}, w_{> K}\big)$ is a product measure, with $\lambda^K$
absolutely continuous with respect to the first factor and with the second
factor equal to $\rho$,
Fubini's theorem gives that for $\rho$-a.e.\ $w$,
for $\lambda^K$-a.e.\ $v \in \R^K$, and for $\lambda$-a.e.\ $s \in \R$, we
have $\bry(\mu)\big((v, w), s\big) = 0$. Since
$(v, s) \mapsto \bry(\mu)\big((v, w), s\big)$ possesses sufficient continuity
properties, we have that for $\rho$-a.e.\ $w$,
for all $v \in \R^K$ and all $s \in \R$,
$\bry(\mu)\big((v, w), s\big) = 0$. 

Let $\epsilon > 0$. Choose $K$ so large that $c \int \|u_{>K}\|_2
\,d\mu_i(u) < \epsilon^2$ for $i = 1,2$, which is possible by Lebesgue's
dominated convergence theorem and the fact that $\mu_i$ has finite first
moment.
Let 
$$
A(\epsilon) := \big\{ (u, w) \in \ell^2(\Z^+) \times \R^\infty
\st |w(u_{>K})| \ge \epsilon \big\}
\,.
$$
Markov's inequality yields that
$$
(\mu_i \times \rho)A(\epsilon) \le 
\epsilon^{-1} \|w(u_{>K})\|_{L^1(\mu_i \times \rho)}
= \epsilon^{-1} c \int \|u_{>K}\|_2 \,d\mu_i(u)
< \epsilon
\,,
$$
where the equality arises from Fubini's theorem.
Therefore, there is some $w$ such that 
denoting $A(w, \epsilon) := \{ u \st |w(u_{>K})|\ge \epsilon \}$,
we have $\bry(\mu)\big((v, w), s\big) = 0$ for all $v \in \R^K$, $s \in \R$
and
$$
\mu_i A(w, \epsilon) < \epsilon
\,.
$$
For such a $w$, we have for all $v, s$ that
$$\eqaln{
\mu_i\big\{u \st \ip u_{\le K}, v \le s - \epsilon\big\} - \epsilon
&<
\mu_i\big\{u \st \ip u_{\le K}, v + w(u_{>K}) \le s \big\}
\cr&<
\mu_i\big\{u \st \ip u_{\le K}, v \le s + \epsilon\big\} + \epsilon
\,.
\cr}$$
The middle quantity is the same for $i = 1$ as for $i = 2$ by choice of $w$.
Therefore, for all $v \in \R^K$ and $s \in \R$,
$$
\mu_1\big[\ip u_{\le K}, v \le s - \epsilon\big] - \epsilon
<
\mu_2\big[\ip u_{\le K}, v \le s + \epsilon\big] + \epsilon
$$
and
$$
\mu_2\big[\ip u_{\le K}, v \le s - \epsilon\big] - \epsilon
<
\mu_1\big[\ip u_{\le K}, v \le s + \epsilon\big] + \epsilon
\,.
$$
Although $K$ depends on $\epsilon$, it follows that for all $L \le K$ and
all $v \in \R^L$, $s \in \R$,  
$$
\mu_1\big[\ip u_{\le L}, v \le s - \epsilon\big] - \epsilon
<
\mu_2\big[\ip u_{\le L}, v \le s + \epsilon\big] + \epsilon
$$
and
$$
\mu_2\big[\ip u_{\le L}, v \le s - \epsilon\big] - \epsilon
<
\mu_1\big[\ip u_{\le L}, v \le s + \epsilon\big] + \epsilon
\,.
$$
Thus, if we fix $L$, the above inequalities hold for all $\epsilon$, which
implies that
$$
\mu_1\big[\ip u_{\le L}, v \le s\big]
=
\mu_2\big[\ip u_{\le L}, v \le s\big]
\,.
$$
This is what we needed to show.
\Qed

Non-separable Hilbert spaces $H$ are of strong negative type iff their
dimension is a cardinal of measure zero. (Whether there exist cardinals not
of measure zero is a subtle question that involves foundational issues; see
Chapter 23 of \ref b.JustWeese/.) To
see this equivalence, note first that if every Borel probability measure on
$H$ is carried by a separable subset, then $H$ has strong negative type by
the preceding theorem. Now a theorem of \ref b.MarczSik/ (or see Theorem 2
of Appendix III in \ref b.Billingsley:conv/) implies that this
separable-carrier condition holds if (and only if) the dimension of $H$ is
a cardinal of measure zero.
Conversely, if the dimension of $H$ is not a cardinal of measure zero, then
let $I$ be an orthonormal basis of $H$. By definition, there exists a
probability measure $\mu$ on the subsets of $I$ that vanishes on
singletons.
Write $I = I_1 \cup I_2$, where $I_1$ and $I_2$ are disjoint and
equinumerous with $I$. Define $\mu_j$ ($j = 1, 2$) on $I_j$ by pushing
forward $\mu$ via a bijection from $I$ to $I_j$. Extend $\mu_j$ to $H$ in
the obvious way (all subsets of $I$ are Borel in $H$ since they are
$G_\delta$-sets).
Then $\mu_1 \ne \mu_2$, yet $D(\mu_1 - \mu_2) = 0$.

\procl c.sqrtnegtype
If $(\Xsp, d)$ is a separable metric space of negative type, then $(\Xsp,
d^{1/2})$ is a metric space of strong negative type.
\endprocl

\proof
Let $\phi : (\Xsp, d^{1/2}) \to H$ be an isometric embedding to a separable
Hilbert space.
Let $\psi : (H, \|\cbuldot\|^{1/2}) \to H'$ be an isometric embedding to
another separable Hilbert space such that $\bry_\psi$ is injective on
$M_1^1(H)$, which exists by \ref t.Hilbertstrong/.
Then $\psi \circ \phi : (\Xsp, d^{1/4}) \to H'$ is an isometric embedding
to a Hilbert space whose barycenter map is injective on $M_1^1(\Xsp,
d^{1/2})$.
\Qed

This means that we can apply a distance covariance test of independence to
any pair of metric spaces of negative type provided we use square roots of
distances in place of distances.
This even has the small advantage that the probability measures in question
need have only finite half-moments.

\procl r.otherpowers
In fact, \ref b.Linde:rudin/ proves that the map $\alpha : \mu \mapsto
a_\mu$ of \ref r.distancetransform/ is injective on $M_1^1(H,
\|\cbuldot\|^r)$ for all $r \in \R^+ \setminus 2\N$. 
It follows that if $(\Xsp, d)$ has negative type, then
$(\Xsp, d^r)$ has strong negative type when $0 < r < 1$.
For let $\phi : (\Xsp, d^{1/2}) \to H$ be an isometric embedding.
By Linde's result, the map 
$$
\mu \mapsto
\left(x \mapsto \int d(x, x')^r \,d\mu(x') = \int
\|\phi(x) - \phi(x')\|^{2r} \,d\mu(x')\right)
$$
is injective.
Since $(\Xsp, d^r)$ has negative type by a theorem of \ref
b.Schoenberg:TAMS/, the claim follows from \ref r.distancetransform/.
\endprocl

\procl c.sqrtsum
If $(\Xsp, d_\Xsp)$ and $(\Ysp, d_\Ysp)$ are metric spaces of negative
type, then $\big(\Xsp \times \Ysp, (d_\Xsp + d_\Ysp)^{1/2}\big)$ is a
metric space of strong negative type.
\endprocl

\proof
It is easy to see that $(\Xsp \times \Ysp, d_\Xsp + d_\Ysp)$ is of
negative type, whence the result follows from \ref c.sqrtnegtype/. 
\Qed

Thus, another way to test independence for metric spaces $(\Xsp, d_\Xsp)$
and $(\Ysp, d_\Ysp)$ of negative type (not necessarily strong) uses not
$\dcov(\theta)$, but $D(\theta - \mu \times \nu)$ with respect to the
metric $(d_\Xsp + d_\Ysp)^{1/2}$ on $\Xsp \times \Ysp$; compare \ref
r.Direct/.
By \ref r.otherpowers/, the same holds for 
$\big(\Xsp \times \Ysp, (d_\Xsp + d_\Ysp)^r\big)$ 
with any $r \in (0, 1)$.

We remark finally that for separable metric spaces of negative type, the
proofs of \ref p.lln/, \ref t.cltgen/, and \ref c.testgen/ are more
straightforward, as they can rely on the strong law of large numbers
and the central limit theorem in Hilbert space.

\medbreak
\noindent {\bf Acknowledgements.}\enspace 
I am grateful to G\'abor Sz\'ekely for asking in what metric spaces a 
theory of distance covariance applies, and for other encouraging
correspondence.
I thank Ori Gurel-Gurevich for a proof of a version of \ref l.tenneg/,
as well as Omer Angel and Assaf Naor for discussions.
I thank Svante Janson for corrections on an earlier draft and a
simplification of the example in \ref r.negstrong/.

\def\noop#1{\relax}
\input \jobname.bbl

\filbreak
\begingroup
\eightpoint\sc
\parindent=0pt\baselineskip=10pt

Department of Mathematics,
831 E. 3rd St.,
Indiana University,
Bloomington, IN 47405-7106
\emailwww{rdlyons@indiana.edu}
{http://mypage.iu.edu/\string~rdlyons/}

\endgroup

\bye